\documentclass[a4paper,12pt]{amsart}

\usepackage{amssymb}
\usepackage{xy}
\usepackage{graphicx}
\xyoption{all}

\newtheorem{theo}{Theorem}[section]

\newtheorem{lemma}[theo]{Lemma}

\newtheorem{conj}[theo]{Conjecture}

\theoremstyle{definition}

\newtheorem{question}[theo]{Question}

\def\deg{\operatorname{deg}\nolimits}
\def\GL{\operatorname{GL}\nolimits}

\def\Disc{\operatorname{Disc}\nolimits}

\def\reg{\operatorname{reg}\nolimits}
\def\sing{\operatorname{sing}\nolimits}
\def\smooth{\operatorname{smooth}\nolimits}

\def\ie{{\em i.e.}}
\def\S{{\mathfrak{S}}}
\def\BP{{\mathbb{P}}}

\def\BC{{\mathbb{C}}}
\def\BN{{\mathbf{N}}}
\def\BQ{{\mathbb{Q}}}

\def\BZ{{\mathbb{Z}}}   
\def\CA{{\mathcal{A}}}
\def\CH{{\mathcal{H}}}
\title{Explicit presentations for exceptional braid groups}
\author{David Bessis, Jean Michel}
\date{}
\address{
DMA, \'Ecole normale sup\'erieure,
45 rue d'Ulm, 75230 Paris cedex 05, France}
\email{david.bessis@ens.fr}
\address{LAMFA, Universit\'e de Picardie-Jules Verne,
33 Rue Saint-Leu, 80039 Amiens Cedex, France\hfill\break\noindent
{\phantom{mm} \sl and} \hfill\break
Institut de Math\'ematiques, Universit\'e Paris VII,
  175 rue du Chevaleret, 75013 Paris, France}
\email{jmichel@math.jussieu.fr}

\thanks{Work on this project was completed during the first author's
stays at Institut Girard Desargues (Lyon) and
Independent University of Moscow.
Computer resources were provided by the Institut de Math\'ematiques
de Jussieu (Paris) -- we thank the computer support team for the
quality of the environment. We thank Enrique Artal Bartolo,
Jorge Carmona Ruber and Bernard Teissier for useful discussions.}
\begin{document}

\begin{abstract}
We give presentations for the braid groups associated
with the complex reflection groups $G_{24}$ and $G_{27}$.
For the cases of $G_{29}$,
$G_{31}$, $G_{33}$ and $G_{34}$, we give (strongly supported)
conjectures.
These presentations were obtained with VKCURVE, a GAP package
implementing Van Kampen's method.
\end{abstract}

\maketitle

\section{Introduction}
To any 
complex reflection group $W\subset\GL(V)$, one may attach a braid
group $B(W)$, defined as the fundamental group of the space
of regular orbits for the action of $W$ on $V$ (\cite{BMR}). 

The ``ordinary'' braid group on $n$ strings, introduced by
Artin (\cite{artin}), corresponds to the
case of the symmetric group $\S_n$, in its monomial reflection representation
in $\GL_n(\BC)$.
More generally, any Coxeter group can be seen as a complex reflection
group, by complexifying the reflection representation. Brieskorn proved
in \cite{brieskorn}
that the corresponding braid group can be described by an \emph{Artin
presentation}, obtained by ``forgetting'' the quadratic relations in
the Coxeter presentation. 

Many geometric properties of Coxeter groups still hold for
arbitrary complex reflection groups. Various authors, including
Coxeter himself, have described ``Coxeter-like'' presentations for
complex reflection groups.
Of course, one would like to have not just a ``Coxeter-like''
presentation for $W$, but also an ``Artin-like'' presentation
for $B(W)$.

The problem can be reduced to the irreducible case.
Irreducible complex reflection groups have been classified
by Shephard and Todd (\cite{ST}): there is an infinite
family $G(de,e,r)$ (which contains the infinite families
of Coxeter groups), plus
$34$ exceptional groups $G_4,\dots,G_{37}$ (among them are the
exceptional Coxeter groups).

Before this note, presentations were known for all but
$6$ exceptional groups (see the tables of \cite{BMR}):
\begin{itemize}
\item The braid group of $G(de,e,r)$ is computed in \cite{BMR}.
The proof makes use of fibration arguments, taking advantage
of the fact that $G(de,e,r)$ is monomial.
\item The first exceptional groups ($G_4$ to $G_{22}$) are
$2$-dimensional. The spaces of regular orbits are complements
of (fairly elementary) complex algebraic curves; the braid
groups have been computed by Bannai (\cite{bannai}),
using Zariski/Van Kampen method.
\item Among the fifteen higher-dimensional exceptional groups,
six are Coxeter groups: Brieskorn's theorem applies to them.
In addition, three more groups happen to have orbit spaces
isomorphic to orbit spaces
of certain Coxeter groups (this was observed by Orlik-Solomon, \cite{OS}).
\item The six remaining groups are $G_{24}$, $G_{27}$,
$G_{29}$, $G_{31}$, $G_{33}$ and $G_{34}$. 
No presentation
for their braid groups are given in \cite{BMR} (except a
conjectural one for $G_{31}$). 
\end{itemize}

In the present note, we describe
presentations for first two of the six remaining groups,
and conjectural presentations for the last four.
The evidence for our conjectures is very strong, and only
a minute step of the proof is missing.

\section{The presentations}

Before listing the individual presentations, it is worth noting that
they share some common features: the number of generators is equal to the
rank of the group (except for $G_{31}$, where an additional generator
is needed); the generators correspond geometrically to 
\emph{generators-of-the-monodromy} (in the sense of \cite{BMR} and
\cite{zariski}) or equivalently \emph{braid reflections}
(this nicer terminology was introduced in \cite{boston});
the relations are positive and homogeneous; by adding quadratic
relations to the presentation, one gets a presentation for the 
reflection group; the product of the
generators, taken in a certain order, has a central power.
Existence of such presentations was proved in \cite{zariski}.
All presentations below satisfy all these properties.

\subsection{The $3$-dimensional group $G_{24}$}

\begin{theo}
The braid group associated with the complex reflection group $G_{24}$
admits the presentation
$$\left< s,t,u \left|
\begin{array}{c}
stst=tsts, susu=usus, tutu=tutu, \\
stustus=tustust=ustustu
\end{array}
\right.
\right> .$$
These relations imply that $(stu)^7$ is central.
\end{theo}

We suggest to represent this presentation by the following
diagram:

$$
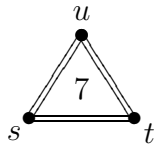
\xy  
(6,11) *++={u} ; (6,8) *=0{\bullet} ;
(-1,-3) *=0{\bullet}  **@{=}; (13,-3) *=0{\bullet} **@{=} ;
(6,8) *=0{\bullet} **@{=};
    (-3,-5) *++={s} ; (15,-5) *++={t} ;
   (6,1) *++={7} \endxy$$

Playing with the above presentation, one may obtain other ones, less
symmetrical but also interesting. E.g., replacing
$t$ by $usts^{-1}u^{-1}$ 
gives (after simplication)

$$<s,t,u | sts=tst,  
          tut=utu,   
          susu=usus, 
          sustustus=ustustust >.$$ 
Also, replacing $t$ by $susts^{-1}u^{-1}s^{-1}$ 
yields
$$<s,t,u | sts=tst,  
          tutu=utut, 
          susu=usus, 
          sutsuts=usutsut >.$$ 

\subsection{The $3$-dimensional group $G_{27}$}

For $G_{27}$, we couldn't find any nice symmetrical presentation, involving
only classical braid relations and cyclic three-terms relations. 

\begin{theo}
The braid group associated with the complex reflection group $G_{27}$
admits the presentations:
\begin{eqnarray*}
< s,t,u & |  & 
stst=tsts,
tut=utu,
sus=usu,
stustustusts
=tstustustust
> \\
< s,t,u & |  &
ststs=tstst,
tut=utu,
sus=usu,
 ststustust=tstustustu > \\
<s,t,u & | &
ststs=tstst,
tutu=utut,
sus=usu,
stustut=
ustustu >.
\end{eqnarray*}
In each of these presentations, the element $(stu)^5$ is central.
\end{theo}

These presentations could be symbolized by the following
diagrams:

$$
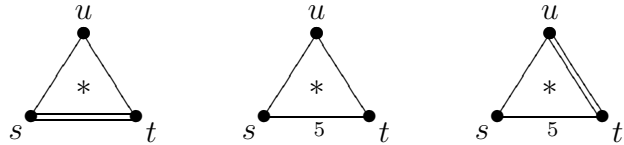
\xy  
(6,11) *++={u} ; (6,8) *=0{\bullet} ;
(-1,-3) *=0{\bullet}  **@{-}; (13,-3) *=0{\bullet} **@{=};
(6,8) *=0{\bullet} **@{-};
    (-3,-5) *++={s} ; (15,-5) *++={t} ;
    (5.9,-5.4) *++={\phantom{.}};
   (6,1) *++={\ast} \endxy \qquad
  \xy  
(6,11) *++={u} ; (6,8) *=0{\bullet} ;
(-1,-3) *=0{\bullet}  **@{-}; (13,-3) *=0{\bullet} **@{-};
(6,8) *=0{\bullet} **@{-};
    (-3,-5) *++={s} ; (15,-5) *++={t} ;
    (5.9,-5.4) *++={\phantom{.}^5};
   (6,1) *++={\ast} \endxy  \qquad
\xy  
(6,11) *++={u} ; (6,8) *=0{\bullet} ;
(-1,-3) *=0{\bullet}  **@{-}; (13,-3) *=0{\bullet} **@{-};
(6,8) *=0{\bullet} **@{=};
    (-3,-5) *++={s} ; (15,-5) *++={t} ;
    (5.9,-5.4) *++={\phantom{.}^5};
   (6,1) *++={\ast} \endxy$$

\subsection{The $4$-dimensional group $G_{29}$}

The presentation for $G_{29}$ given in \cite{BMR} was not conjectured
to give (by forgetting the quadratic relations) a presentation
for the braid group. Surprisingly, our computations happened to
give precisely this presentation.

\begin{conj}
The braid group associated with the complex reflection group $G_{29}$
admits the presentation
$$\left< s,t,u,v \left|
\begin{array}{c}
sts= tst, tut=utu, uvu=vuv, tvtv=vtvt,  \\
su=us, sv=vs, utvutv=tvutvu
\end{array}
\right.
\right>.$$
These relations imply that $(stuv)^5$ is central.
\end{conj}

Brou\'e-Malle-Rouquier used the following diagram to symbolize
this presentation:

$$
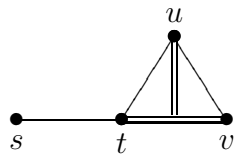
\xy  
(6,11) *++={u} ; (6,8) *=0{\bullet} ;
(-1,-3) *=0{\bullet}  **@{-}; (13,-3) *=0{\bullet} **@{=} ;
(6,8) *=0{\bullet} **@{-};
(-1,-3) *=0{\bullet}; (-15,-3) *=0{\bullet} **@{-};
    (-1,-6) *++={t} ; (13,-6) *++={v} ;
    (-15,-6) *++={s} ;
 (6,8) *=0{\bullet} ; (6,-4) *++={\phantom{.}} **@{=}  
\endxy$$

\subsection{The $4$-dimensional group $G_{31}$}

The following conjecture ``confirms'' the conjectural presentation
from \cite{BMR} -- but this time there is computational evidence behind.

\begin{conj}
The braid group associated with the complex reflection group $G_{31}$
admits the presentation
$$\left< s,t,u,v,w \left|
\begin{array}{c}
sts= tst, tut=utu, uvu=vuv, vwv=wvw \\
sv=vs, tv=vt, tw=wt \\
suw=uws=wsu
\end{array}
\right.
\right>.$$
These relations imply that $(stuwv)^6$ is central.
\end{conj}

The corresponding
Brou\'e-Malle-Rouquier diagram is:

$$
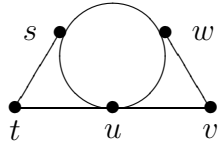
\xy  
(2,7) *++={s} ; (25,7) *++={w} ; (6,7) *=0{\bullet} ;
(0,-3) *=0{\bullet}  **@{-}; (13,-3) *=0{\bullet} **@{-} ;
 (26,-3) *=0{\bullet}  **@{-}; (20,7) *=0{\bullet} **@{-} ;
   (0,-6) *++={t} ; (13,-6) *++={u} ; (26,-6) *++={v} ;
   (13,3.8) *++={\phantom{sdmds}} *\frm{o} \endxy$$

{\flushleft 
\bf Remark.} Since our generators are braid reflections, they
map to generating reflections in the reflection group. It is well-known
that, even though it is $4$-dimensional,
$G_{31}$ cannot be generated by less than $5$ reflections.

\subsection{The $5$-dimensional group $G_{33}$}
The relations in the presentation below do not coincide with the
homogeneous part of the Brou\'e-Malle-Rouquier presentation of
$G_{33}$. However, the relations involving $t,u,w$ coincide with the
Brou\'e-Malle-Rouquier relations for the braid group of $G(3,3,3)$
(the similar remarks also apply to $G_{34}$).
\begin{conj}
The braid group associated with the complex reflection group $G_{33}$
admits the presentation
$$\left< s,t,u,v,w \left|
\begin{array}{c}
sts=tst, tut=utu, uvu=vuv, wtw=twt, wuw=uwu \\
su=us, sv=vs, tv=vt, ws=sw, wv=vw \\
tuwtuw=uwtuwt=wtuwtu
\end{array}
\right.
\right>.$$
These relations imply that $(stuvw)^9$ is central.
\end{conj}

(the relation $uwtuwt=wtuwtu$ is redundant).

We suggest to represent this presentation by the following
diagram:

$$
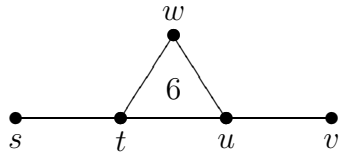
\xy  
(6,11) *++={w} ; (6,8) *=0{\bullet} ;
(-1,-3) *=0{\bullet}  **@{-}; (13,-3) *=0{\bullet} **@{-} ;
(6,8) *=0{\bullet} **@{-};
(-1,-3) *=0{\bullet}; (-15,-3) *=0{\bullet} **@{-};
(13,-3) *=0{\bullet}; (27,-3) *=0{\bullet} **@{-};
    (-1,-6) *++={t} ; (13,-6) *++={u} ;
    (-15,-6) *++={s} ; (27,-6) *++={v} ;
   (6,1) *++={6} \endxy$$

Following \cite{BMR} where a second diagram for $G_{33}$ is given
(to account for some parabolic subgroups missing in their first
diagram), it is not difficult to obtain the equivalent presentation
$<s,t,u,v,w | vt=tv,uv=vu,tu=ut,wu=uw,wsw=sws,sus=usu,svs=vsv,
sts=tst,vwv=wvw,twt=wtw,twvstw=wvstwv>$, which 
contains a parabolic subdiagram of type $D_4$.
(A similar diagram may be derived from the conjectural presentation
for $B(G_{34})$ given below).
   
\subsection{The $6$-dimensional group $G_{34}$}
\begin{conj}
The braid group associated with the complex reflection group $G_{34}$
admits the presentation
$$\left< s,t,u,v,w,x \left|
\begin{array}{c}
\text{relations of $G_{33}$ +} \\
xvx=vxv,xs=sx,xt=tx,xv=vx,xw=wx
\end{array}
\right.
\right>.$$
These relations imply that $(stuvwx)^7$ is central.
\end{conj}

We suggest to represent this presentation by the following
diagram:

$$
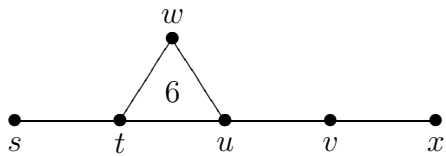
\xy  
(6,11) *++={w} ; (6,8) *=0{\bullet} ;
(-1,-3) *=0{\bullet}  **@{-}; (13,-3) *=0{\bullet} **@{-} ;
(6,8) *=0{\bullet} **@{-};
(-1,-3) *=0{\bullet}; (-15,-3) *=0{\bullet} **@{-};
(13,-3) *=0{\bullet}; (27,-3) *=0{\bullet} **@{-};
(41,-3) *=0{\bullet} **@{-};
    (-1,-6) *++={t} ; (13,-6) *++={u} ;
    (-15,-6) *++={s} ; (27,-6) *++={v} ;
    (41,-6)  *++={x} ; 
   (6,1) *++={6} \endxy$$

\section{Definitions and preliminary work}

Our strategy of proof is, basically, brute force.
Let $V$ be a $\BC$-vector space of dimension $r$, and let
$W\subset\GL(V)$ be a complex reflection group.
The algebra $\BC[V]^W$ of invariant polynomial functions is isomorphic
to a polynomial algebra (\cite{ST}); let $f_1,\dots,f_r$ be homogeneous
polynomials such that $\BC[V]^W=\BC[f_1,\dots,f_r]$.

Let $\CA$ be the set of all reflecting hyperplanes.
For each $H\in\CA$, the pointwise stabilizer $W_H$ of $H$ in $W$ is
a cyclic subgroup of order $e_H$; choose $l_H$ a linear form with kernel
$H$. Let $V^{\reg} :=V-\bigcup_{H\in \CA}H$.
The \emph{regular orbits space} is $V^{\reg}/W$.
We have $\prod_{H\in\CA}l_H^{e_H}\in \BC[V]^W$, so there is a unique
polynomial $\Delta\in\BC[X_1,\dots,X_r]$ such that
$\prod_{H\in\CA}l_H^{e_H} = \Delta(f_1,\dots,f_r)$.
We call $\Delta$ the \emph{discriminant} of $W$ (with respect to
$f_1,\dots,f_r$).
Clearly, $V^{\reg}/W$ is isomorphic, as an algebraic variety, to
the complement of the hypersurface $\CH$
defined in $\BC^r$ by the equation
$\Delta=0$.

There is a general method, though not always practically
tractable, to compute
the fundamental group of such a space.
First, choose a $2$-plane $P$ such that the embedding
$P\cap(\BC^r-\CH) \hookrightarrow \BC^r-\CH$ is a $\pi_1$-isomorphism
(by a Zariski theorem,
this should hold for a \emph{generic} choice of $P$ -- how exactly
this choice can be made is a difficult issue, which we will discuss
later on).
Then use the Zariski/Van Kampen method to compute the fundamental
group of $P\cap(\BC^r-\CH)$.
The computations involved in the second step are far beyond
human capabilities (or at least beyond \emph{our} capabilities),
especially if one wants to avoid imprecise arguments.
Therefore we designed a software package, VKCURVE (\cite{VK}), to carry 
them by computer.

\subsection{General remarks about the implementation}
Our computations are performed using the 
computer
algebra software GAP3, which is designed to handle 
cyclotomic numbers, matrices over these numbers, permutations,
presentations, and all sorts of algebraic objects and algorithms
involving exact computations. The source of its
mathematically advanced functions is public (and in a rather intelligible
language) and any user is free to check their validity.

Our package VKCURVE builds on the older package CHEVIE, which implements
(among other)
complex reflection groups, Coxeter groups and Artin groups.

\subsection{Computing the discriminant}
For each of  the six groups, the discriminant can  be recovered from the
data given in Appendix B of \cite{OT}, where Orlik and Terao explain how
to construct  the matrix $M \in  M_r(\BC[X_1,\dots,X_r])$ of logarithmic
vector  fields (aka. basic derivations)
for  the  quotient singularity  (called the  discriminant
matrix  in  \cite[6.67]{OT}).  The  polynomial $\Delta$  is  simply  the
determinant of $M$.

To prevent typos, we actually re-checked all needed computations.

We   summarize    their   method.   Let   $d_i=\deg    f_i$,   and   let
$d^*_1,\ldots,d^*_j$ be the codegrees of  $W$. We assume that the degrees are
ordered in  increasing order (but  we {\it do  not} assume the  same for
codegrees). The matrix $M$ is  an $r\times r$-matrix whose $(i,j)$ entry
is an homogeneous invariant polynomial of degree $d_i+d^*_j$.

The six groups have the property that $d_1<d_2$ so $f_1$ is unique up to
a  scalar, and  if  $H=(\partial_j\partial_i f_1)_{ij}$  is the  Hessian
matrix  of $f_1$,  then $\det  H$  may be  chosen  as one  of the  basic
invariants $f_k$  (which we assume). Then,  if $J=(\partial_j f_i)_{ij}$
is  the Jacobian  matrix  of the  $f_i$, we  have  the following  matrix
equation (\cite[(1)  p.280]{OT}): $$ M=(d_1-1)JH^{-1}{}^tJ  C\eqno (1)$$
where $C$  is a  matrix of homogeneous  invariant polynomials  such that
$\deg C_{ij}=d_1+d^*_j-d_i$. Orlik  and Terao note that  there exists an
ordering of the  $d^*_j$ such that $C$ is the  identity matrix, except
for some line $q$ where $C_{iq}=0$ for $i<q$, $C_{qq}=f_k$, and $C_{iq}$
is  a polynomial  in  $f_1,\ldots,f_{k-1},f_{k+1},\ldots,f_r$ for  $i>q$
(the  degrees  of the  entries  of  $C$ determine the
ordering).

Equation (1) is used first to  determine $C$, and then to determine $M$.
It  may  be used  to  determine  $C$  since  it implies  the  polynomial
congruence $0\equiv JH'  {}^tJC \pmod{f_k}$, where $H'$  is the cofactor
matrix of  $H$; each non-zero  entry of $C$  is a linear  combination of
(known from  their degree)  monomials in the  basic invariants,  and the
above polynomial congruence is  sufficient to determine the coefficients
of the linear combination.

{\small {\bf \flushleft Example.}
Sufficient data to construct the matrix of basic derivations for $G_{24}$
is given on p. 284 in \cite{OT}.
Note however that the formula given p. 264 in \cite{OT} for its determinant
contains a typo. The correct formula is
\begin{eqnarray*}
\Delta_{24} & = & -2048 x^9 y + 22016 x^6y^3 - 60032 x^3y^5 + 1728 y^7
- 256 x^7z  \\ & &  + 1088 x^4y^2z  + 1008 xy^4z^2 - 88 x^2yz^2 + z^3.
\end{eqnarray*}
To check that such a formula is correct, it suffices to substitute
the invariants: the result should be the product of the
square of the linear forms defining the reflecting hyperplanes.}


In the Appendix, we list basic derivations for all examples (except
$G_{34}$, for which the matrix is too large to be printed...)

\section{Choosing the $2$-plane}

\subsection{A general strategy}
An explicit \emph{genericity} criterion is given in \cite{dimca}, Ch. 4,
Theorem 1.17: it suffices that $P$ is transverse to all the strata of
a Whitney stratification of the hypersurface.
The theorem applies to a projective context.
We replace $\Delta$ by a homogeneous polynomial
$\tilde{\Delta}\in\BC[X_0,\dots,X_r]$.
The equation $\tilde{\Delta}=0$ defines a projective
hypersurface $\tilde{\CH}$; we are interested in the complement
$\BC^r-\CH =  \BC\BP^r-\tilde{\CH}\cup \BC\BP^{r-1}$.

First, we stratify $\BC^r$ as follows:
for all $k\in \{0,\dots,r\}$, set $E_k$ to be the locus
where the matrix $M$ has rank $k$.
This stratification is the quotient modulo $W$ of the stratification
of $V$ by the intersection lattice of $\CA$, hence is a Whitney 
stratification. Moreover, the tangent space of the stratum at a given
point is spanned by the columns of $M$.
With the explicit knowledge of this matrix, there is no major difficulty in
checking transversality of a given $2$-plane.

{\small 
{\flushleft \bf Example.}
For $G_{31}$, one may check that the transversality at infinity
is statisfied by
the $2$-plane of the equations
\begin{eqnarray*}
z & = & y \\
t & = & 1+x
\end{eqnarray*}
The affine tranversality condition for this $2$-plane
is that, for each value of $x$ and $y$,
the following matrix has rank $4$ (the matrix of basic derivations
for $G_{31}$ is given in the Appendix):

$${ \begin{pmatrix}
\scriptstyle 8x & \scriptstyle 12y+12xy &
\scriptstyle 20y+\frac{1}{135}xy & 
\scriptstyle 24+24x-\frac{1}{135}y^2  & 
\scriptstyle 1 & 
\scriptstyle 0 \\
\scriptstyle  12y  & 
 \scriptstyle 18x^2-97200y^2+18x^3 & 
 \scriptstyle -36-36x+\frac{1}{90}x^3 & \scriptstyle -42xy-\frac{1}{90}x^2y 
 & \scriptstyle 0 & \scriptstyle 1\\
\scriptstyle   20y &  \scriptstyle 
-36-72x+60xy^2-36x^2 & \scriptstyle -\frac{1}{270}x-\frac{1}{270}x^2-
\frac{1}{54}y^2 &
    \scriptstyle   \frac{1}{270}y+\frac{1}{270}xy+\frac{1}{54}x^2y  &
    \scriptstyle 0 & \scriptstyle 1 \\
\scriptstyle  24+24x & \scriptstyle -42xy-42x^2y-60y^3
& \scriptstyle -\frac{1}{54}y-\frac{1}{54}xy-\frac{1}{135}x^2y &
   \scriptstyle    \frac{1}{135}xy^2+\frac{1}{54}x^2+20y^2+\frac{1}{54}x^3 
   & \scriptstyle 1 & \scriptstyle 0
\end{pmatrix}}$$
where the first four columns generate the tangent vector to the local
stratum of the discriminant and the last two columns generate
the tangent vector space to the $2$-plane.}

To apply \cite[Ch. 4, Theorem 1.17]{dimca}, we also need a stratification
of the hyperplane at infinity $\BC\BP^{r-1}$.
Let $\CH_{\infty}:=\tilde{\CH}\cap \BC\BP^{r-1}$.
We view $\CH_{\infty}$ as an algebraic hypersurface in $\BC\BP^{r-1}$, defined
by the equation $\sqrt{\Delta_{\infty}}=0$, where $\Delta_{\infty}$
is the homogeneous part of highest degree of $\Delta$, and 
$\sqrt{\Delta_{\infty}}$ is a reduced version of $\Delta_{\infty}$.
We set $N_{r-1}:=\BC\BP^{r-1} - \CH_{\infty}$,
$N_{r-3}:=(\CH_{\infty})_{\text{sing}} \cup \overline{M_{r-2}}$
and $N_{r-2}:=\CH_{\infty} - N_{r-3}$.

Together, the $M_i$'s and the $N_i$'s form a stratification (without
border condition), with incidence
diagram:

$$
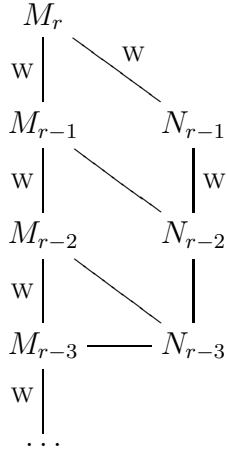
\xymatrix{
M_r \ar@{-}[d]_{\text{W}} \ar@{-}[dr]^{\text{W}} \\
M_{r-1} \ar@{-}[d]_{\text{W}} \ar@{-}[dr] & N_{r-1} \ar@{-}[d]^{\text{W}} \\
M_{r-2} \ar@{-}[d]_{\text{W}} \ar@{-}[dr] & N_{r-2} \ar@{-}[d] \\
M_{r-3} \ar@{-}[d]_{\text{W}} \ar@{-}[r] & N_{r-3}  \\
\dots }$$
We mark W where we know that the incidence satisfies Whitney's conditions.
We have already explained why the first column is a Whitney
stratification. It is trivial that $M_r$ is Whitney over $N_{r-1}$ 
and that $N_{r-1}$ is Whitney over $N_{r-2}$.
By splitting $N_{r-3}$ into smaller strata, we may ensure that
everything below $N_{r-2}$ and $M_{r-2}$ is Whitney (see for example
the construction explained at the beginning of \cite{GWPL}).

\begin{question}
\label{laquestion}
Does $M_{r-1}$ satisfy Whitney's conditions over $N_{r-2}$?
\end{question}

Note that, since $\Delta_{\infty}$ is not reduced, the points
of $N_{r-2}$ are not smooth in $\CH$, so the answer is not that trivial.
It is a pity that no software is available to answer such a question, on
specific examples with explicit equations.

{\small
{\bf \flushleft Example.} For $G_{31}$, we represent points of
$\BC\BP^5$ by $5$-tuples $(h,x,y,z,t)$, with either $h=1$
(affine portion) or $h=0$ (space at infinity).
The strata $M_i$ have explicit
equations, using the matrix given in the Appendix.
The affine hypersurface $\CH$ is given by $h=1$ and
$\Delta(x,y,z,t)=0$, where $\Delta$ is the determinant
of the relevant matrix from the Appendix.
We have 
$\Delta_{\infty}=
-\frac{4}{27}x^7z^2t-\frac{8}{81}x^6yz^3$, thus
a reduced equation for $\CH_{\infty}$ is
$xz(3xt+2yz)$ (and $h=0$). 
One may prove (by means of
Gr\"obner basis) that if a sequence $(1,x_m,y_m,z_m,t_m)_{m\in\BN}$ of
points in
$M_2\cup M_1\cup M_0$ converges to $(0,\overline{x},\overline{y},
\overline{z},\overline{t})$, then either $(\overline{x},\overline{y})=(0,0)$,
or $(\overline{x},\overline{z})=(0,0)$, 
or $(\overline{t},\overline{z})=(0,0)$.
This locus actually coincides with $(\CH_{\infty})_{\sing}$,
thus $N_{r-3}$ is the complement in $\CH_{\infty}$ of this locus
(this explains why the particular
$2$-plane given earlier avoids $N_{r-3}$: the points at infinity of the
$2$-plane have the form $(0,x,y,y,x)$, where either $x\neq 0$ or $y\neq 0$).
Question \ref{laquestion} specializes to: is $\CH_{\smooth}$ Whitney
over $(\CH_{\infty})_{\smooth}$?

}

We may now explain what we have checked, and what is missing to turn our
conjectures into theorems:
\begin{itemize}
\item For all six examples, our presentations
were obtained by applying Van Kampen's method to the algebraic curves
obtained with particular $2$-planes.
\item For all six examples, we have checked that the $2$-planes are
transversal to the affine strata $M_0,\dots,M_r$.
\item For all examples but $G_{34}$,
we have computed (by means of Gr\"obner basis) equations for
$N_{r-1}$, $N_{r-2}$ and $N_{r-3}$, and checked that our $2$-planes are also 
transversal to these strata. Transversality implies that the $2$-planes
do not intersect $N_{r-3}$, and therefore remain transversal to the
Whitney refinement of $N_{r-3}$.
Therefore, if Question \ref{laquestion} had
a positive answer, our conjectures would be theorems (except for $G_{34}$).
\item Note that it is easy to check that our $2$-planes give
generators of the fundamental group, and any homotopy
in the $2$-plane is a homotopy in $\BC^r$. Therefore, we know
for sure that there are presentations for the braid groups
obtained by adding relations to our conjectural presentations. On the
other hand, we have checked that
adding quadratic relations to our conjectural presentations yields
actual presentations for the complex reflection group. Any missed
relation should be trivial in this quotient.
\end{itemize}

\subsection{A strategy for $3$-dimensional groups}
Another approach, more algebraic, can be used to find good $2$-planes.
Although we may start the discussion with any of our examples,
it will be conclusive only for $3$-dimensional groups.
We work with the setting and notations from \cite[Section 2.2]{zariski}:
we have $\Delta\in\BC[X_1,\dots,X_r]$ ($\Delta$ plays the part of the
polynomial $P$ of \emph{loc. cit.}).
We distinguish the variable $X:=X_r$,
we choose a generic (in the sense of \emph{loc. cit.})
line $L$ of direction $X$. Viewed as a polynomial in only the variable $X$
(with coefficients involving the other variables), $P$ has a
discriminant $\Disc(P_X)$. Let
$E:=\{v\in\BC^r |P(v) \neq 0, \Disc(P_X)(v)\neq 0\}$.
We denote by $p$ the projection $(x_1,\dots,x_r) \mapsto (x_1,\dots,x_{r-1})$.
Let $\overline{E}:=p(E)$.
The map $p$ induces a fibration $E \twoheadrightarrow \overline{E}$,
whose exact sequence ends as follows:
$$\xymatrix@1{\cdots \ar[r] & \pi_2(\overline{E},y_0) \ar[r]  &
\pi_1(L-L\cap\CH,x_0) \ar^{\phantom{mu}\iota_*}[r]
& \pi_1(E,(x_0,y_0)) \ar^{\phantom{mu}p_*}[r] &
\pi_1(\overline{E},y_0) \ar[r] & 1}$$
In our setting, $\Delta$
is monic in $X$ (since $d_r$ is regular, it follows from
\cite[Lemma 1.6]{zariski}). It is then easy to construct a section
$s:\pi_1(\overline{E},y_0) \rightarrow \pi_1(E,(x_0,y_0))$ of $p_*$.

The basespace $\overline{E}$ is the complement in $\BC^{r-1}$ of the
hypersurface of equation $\Disc(\Delta_X)=0$. In our setting, 
$\Disc(\Delta_X)$ is a weighted homogeneous polynomial.

When $r=3$, this implies that $\pi_2(\overline{E},y_0)=0$
(complements of weighted homogenenous curves are $K(\pi,1)$).
We then have a semi-direct product structure
$$\pi_1(E,(x_0,y_0)) \simeq
\pi_1(L-L\cap\CH,x_0) \rtimes \pi_1(\overline{E},y_0).$$
To obtain a presentation for $\pi(\BC^r - \CH)$, one
starts with a presentation for $\pi_1(E,(x_0,y_0))$, and adds
relations forcing elements of $\pi_1(\overline{E},y_0)$ to become trivial.
It is an easy exercise to check, in this setting, that any $2$-plane $P$
satisfying:
\begin{itemize}
\item[(i)] the line $L$ is contained in $P$ and,
\item[(ii)] the image line $p(P)$ is such that
$p(P) \cap \overline{E} \hookrightarrow \overline{E}$ is $\pi_1$-surjective,
\end{itemize}
is good for our purposes. In our examples, it is easy to construct such
planes, since
$\Disc(\Delta_X)$ is monic in one the remaining variables.
This is how we obtained, for $G_{24}$ and $G_{27}$,
theorems rather than conjectures.

Note that, for other groups, all assumptions used here (including
the monicity of $\Disc(\Delta_X)$) remain valid, except that
we do not know whether $\pi_2(\overline{E},y_0)=0$. Instead
of answering Question \ref{laquestion}, checking that
$\pi_2(\overline{E},y_0)=0$ would turn our conjectures into theorems.

\section{The package VKCURVE}
Once a $2$-plane $P$ has been chosen, it is enough to feed
VKCURVE
with the equation of the curve $P\cap \CH$ to obtain
a presentation of $\pi_1(P-(P\cap\CH))$.

{\small
{\bf \flushleft Example.} For $G_{31}$, when computing 
the determinant of $M_{31}$ and evaluating at $z=y$ and $t=1+x$,
we obtain the following equation for $P\cap \CH$:
\begin{align*}
\begin{split}
\scriptstyle \Delta'_{31} \; = \;
& \scriptstyle 
746496+3732480x-3111936xy^2-\frac{93281756}{27}xy^4+\frac{58341596}{27}xy^6
+7464960x^2-
384y^2-
9334272x^2y^2 \\ 
& \scriptstyle 
+\frac{17556484}{27}x^2y^4+\frac{43196}{27}x^2y^6+7464576x^3
-\frac{756138248}{81}x^3y^2+\frac{192792964}{81}x^3y^4+\frac{16}{81}x^3y^6
+3730944x^4 \\
& \scriptstyle 
-\frac{139967996}{81}y^4-\frac{84021416}{27}x^4y^2+
\frac{82088}{27}x^4y^4+744192x^5+\frac{43192}{27}x^5y^2-
\frac{1720}{27}x^5y^4 
-\frac{124412}{81}x^6 \\
& \scriptstyle 
+777600800y^6+\frac{95896}{81}x^6y^2
-\frac{8}{81}x^6y^4-\frac{10364}{27}x^7-\frac{4}{27}x^7y^2+\frac{4}{27}x^8
-\frac{8}{81}y^8-\frac{4}{27}x^8y^2+\frac{4}{81}x^9
\end{split}
\end{align*}
On a 3 GHz Pentium IV, VKCURVE needs about one hour to deal with this example.}

Writting VKCURVE was of course the most difficult part of our work.
This software accepts
as input any quadratfrei polynomial in $\BQ[i][X,Y]$ and computes
a presentation for the fundamental
group of the complement of the corresponding complex algebraic
curve. The program does not use floating point computations
(even when computing monodromy braids); therefore
there is no issue of numerical accuracy 
and the result is ``certified'' to be correct (provided that our
implementation does not contain mathematical errors...)

The remainder of this section is an overview of the algorithms
used in VKCURVE. We rely on the version of Van Kampen's 
method exposed in \cite[Procedure 4]{variations}, where it is
decomposed into four steps.

\subsection{Implementing steps 1 and 2}
Starting with our polynomial $P\in\BQ[i][X,Y]$, we view it as
a one variable polynomial in $\BQ[i][Y][X]$ and compute
its discriminant $\Delta\in\BQ[i][Y]$.
The discriminant $\Delta$ may not be reduced; to compute
approximations $\tilde{y}_1,\dots,\tilde{y}_r\in\BQ[i]$
of its complex roots $y_1,\dots,y_r$, we apply Newton's
method to the reduced polynomial $\Delta_0$ obtained by
dividing $\Delta$ by the resultant of $\Delta$ and $\Delta'$.
As Hubbard, Schleicher and Sutherland proved in their beautiful
article \cite{hss}, Newton's method can be made into a failsafe
algorithm producing arbitrarily good approximations of $y_1,\dots,y_r$.

Since we will re-use them later, we recall a
few trivialities about complex polynomials.
Let $P\in \BC[Z]$.
Let $\alpha_1,\dots,\alpha_n$ be the complex roots of $P$.
Let $z\in \BC$. If $P'(z)\neq 0$, we set
$N_P(z) := z - \frac{P(z)}{P'(z)}$.
Considering the first order approximation of $P$ around $z$, we
expect $P(N_P(z))$ to be close to $0$. Newton's method consists of
starting with $z_0\in\BC$ (chosen randomly, or smartly as in
\cite{hss}) and to construct iteratively $z_{m+1}:=N_P(z_m)$,
hoping that $(z_m)$ will
converge towards a root of $P$ -- which indeed happens for
``many'' choices of $z_0$. How may we decide that a given
$z_n$ is a ``good enough'' approximation?

\begin{lemma}
\label{trivial}
Assume $P$ has $n$ distinct roots $\alpha_1,\dots,\alpha_n$. 
Let $z\in \BC$, with $P'(z)\neq 0$.
Then there exists $\alpha\in\{\alpha_1,\dots,\alpha_n\}$ such that
$|z - \alpha| \leq n \left| \frac{P(z)}{P'(z)} \right|$.
\end{lemma}

\begin{proof}
If $P(z)=0$, the result is trivial. Otherwise, we have
$\frac{P'(z)}{P(z)} = \sum_{i=1}^n \frac{1}{z-\alpha_i}$.
Choose $i$ such that for all $j$, $|z-\alpha_i|\leq |z-\alpha_j|$.
By triangular inequality,
$\frac{1}{|z-\alpha_i|} \geq \left| \frac{P'(z)}{P(z)} \right|
- \sum_{j\neq i} \frac{1}{|z-\alpha_j|}
\geq  \left| \frac{P'(z)}{P(z)} \right| -
(n-1) \frac{1}{|z-\alpha_i|}$. The result follows.
\end{proof}

Although elementary, this lemma provides a very inexpensive
(in terms of computational time) test for deciding whether a tentative
list $\tilde{\alpha}_1,\dots,\tilde{\alpha}_n$ of complex numbers
``separates'' the roots (\ie, whether there exists
$\varepsilon_1,\dots,\varepsilon_n$ such the disks
$D(\tilde{\alpha_i},\varepsilon_i)$
do not overlap and each of them contains a root of $P$).

Instead of working with the exact Newton's method, we use a truncated
version, where $N_P(z)$ is replaced by an approximate
$(a+ib)10^k$, where $a,b\in \BZ$, and $k$ is an integer
slightly smaller than $\log_{10} \left| \frac{P'(z)}{P(z)} \right|$.
This is to avoid the very fast increase of the denominators, when
the exact method is carried out in $\BQ[i]$: the complexity
of the exact method is very good from the ``abstract'' viewpoint
(the number of iterations), but in practice really bad
(each individual iteration involves costly operations on very big integers).
Of course, our modification does not make the method less rigorous, since
the test can be performed exactly. The main difference between our
implementation and floating point is that $k$ is modified dynamically
and has no pre-assigned bound.

\smallskip

Once separating approximates
$\tilde{y}_1,\dots,\tilde{y}_r\in\BQ[i]$
of the roots of $\Delta$ have been obtained,
Step 2 of
\cite[Procedure 4]{variations} is performed as follows:
first, we construct the Voronoi cells around
$\tilde{y}_1,\dots,\tilde{y}_r$; then, concatenating some
of the affine
segments bounding the Voronoi cells, we construct, for each $i$, a loop
$\gamma_i$ representing a meridien around
$\tilde{y}_i$; it is easy to make sure
that we recover a meridien around the actual ${y}_i$.

\subsection{Step 3: computing monodromy braids}
\cite[Procedure 12]{variations} decomposes Step 3 into smaller
steps a--e. Only Substep $a$ is not a straightforward algebraic
manipulation -- and most of the computational time is spent there.
The problem is as follows: let $[y_0,y_1]$ be one
of the affine segments involved in the $\gamma_i$'s.
For $t\in [0,1]$, denote by $P_t\in \BQ[i][X]$ the polynomial obtained by
evaluating $P$ at $Y=(1-t)y_0+ty_1$.
We want to compute the word in Artin generators corresponding
to the real projection of the braid obtained by tracking the
roots of $P_t$ when $t$ runs over $[0,1]$.

As we have seen above, we may find $x_1,\dots,x_n\in \BQ[i]$
separating the roots of $P_0$. Concretely, using Lemma
\ref{trivial}, we iterate a truncated Newton method until,
when we set $\varepsilon_i:=\inf_{j\neq i} \frac{|x_i-x_j|}{2}$
(this is a simple way, though not optimal, to ensure that
$\forall i, j, |x_i-x_j| > \varepsilon_i + \varepsilon_j$),
we have
$$\forall i, \left| \frac{P_0(x_i)}{P'_0(x_i)}\right| < 
\frac{\varepsilon_i}{n}.$$
For each $i$, consider the polynomial
$$Q_i:= \varepsilon_i^2 |P'_t(x_i)|^2 - n^2|P_t(x_i)|^2 \in \BQ[t].$$
By assumption, we have $\forall i, Q_i(0) > 0$.
Whenever $t_0\in [0,1]\cap \BQ$ is such that $\forall t\in [0,t_0], \forall i,Q_i(t)>0$,
we know that, for $t\in [0,t_0]$,
the strings of the monodromy braids will be in the cylinders
of radius $\varepsilon_i$ around the $x_i$'s.
This fragment of the monodromy braid can be replaced by the
constant braid with strings fixed at the positions given by the $x_i$'s.
Set $y_0':=(1-t_0)y_0+y_1$,
$x_i':=N_{P_{t_0}}(x_i)$. Though the $x_i$'s already separate the
roots of $P_{t_0}$, the $x_i'$'s shoud be ``better'' approximates.
We compute new radii $\varepsilon_i'$ separating the $x_i$'s and
iterate, studying now the monodromy braid over
$[y_0',y_1]$, with initial approximates $x_1',\dots,x_n'$.
Eventually, we hope that after some number of iterations,
$t_0=1$ will suit.

The main difficulty is to find an actual $t_0$
such $\forall t\in [0,t_0], \forall i,Q_i(t)>0$.
One the one hand, we want it to be as large as possible, to avoid
unnecessary iterations; on the other hand, computing the largest
theoretical value for $t_0$, for example using Sturm sequences,
is very costly. Finding a good balance is a delicate art.
The curious reader may have a look at the source of the
VKCURVE function {\tt FollowMonodromy}, where a very naive
method is used, together with careful coding and adaptative
heuristics (note that, in {\tt FollowMonodromy}, one actually computes
a distinct $t_0$ for each individual string -- the above description
is simplified for the sake of clarity).

\subsection{Step 4: writing and simplifying the presentation}
Working with GAP, it is then straightforward to write
a presentation. However, this presentation
is much more complicated than desirable.
Since no ``normal form'' theory exists for arbitrary
group presentations, it is not clear how one can simplify
it and obtain
one of our ``pretty'' presentations. Fortunately, some natural heuristics
(typically, replace a generator by its conjugate by another generator,
try to simplify, iterate in the regions of the tree of all possibilities
where the total length
of the presentation tends to decrease) happen to be quite
effective in dealing with the
(highly redundant) presentations obtained with Van Kampen's method.
Playing with these (non-deterministic) heuristics, which are part of
VKCURVE, we obtained quite easily a few ``simple'' presentations.
At this point, in the absence of a general combinatorial theory of
generalized braid groups, there
is some arbitrary in deciding which one should be retained;
in most cases though, one of them clearly emerged as being the
``prettiest''.

\appendix

\section{Explicit matrices of basic derivations}

\def\Hessian{\operatorname{Hessian}}
\def\bord{\operatorname{bord}}

\subsection*{The group $G_{24}$}
With Klein's matrices (as in \cite{OT})
the first invariant is $f_1=x^3y+zy^3+xz^3$.
The others are $f_2=\det(\Hessian(f_1))/108$ and $\bord(f_1,f_2)/36$.

Basic derivations:
$$\begin{pmatrix}
4x&6y^{2}&14z-36x^{2}y\\
6y&-z&128xy^{2}-7x^{4}\\
14z&128xy^{3}-6x^{2}z-7x^{4}y&287xyz-35x^{3}y^{2}-294y^{4}+7x^{6}
\end{pmatrix}$$

\subsection*{The group $G_{27}$}
With Wiman's matrices (as in \cite{OT}),
the first invariant is
$f_1=-135xyz^{4}-45x^{2}y^{2}z^{2}+10x^{3}y^{3}+9x^{5}z+9y^{5}z+27z^{6}$.
The others are $f_2=\det(\Hessian(f_1))/6750$ and $\bord(f_1,f_2)/5400$.

Basic derivations:
{\small
$$\begin{pmatrix}
6x&12y^{2}&30z+234x^{3}y\\
12y&-4z&-34xz-1362x^{2}y^{2}+156y^{3}+900x^{4}y-270x^{6}\\
30z&
\stackrel
{\displaystyle -34xyz-1362x^{2}y^{3}+78x^{3}z+}
{\displaystyle 156y^{4}+900x^{4}y^{2}-270x^{6}y}
&
\stackrel
{\displaystyle 4836xy^{4}+330y^{2}z-3349x^{2}yz-17727x^{3}y^{3}+}
{\displaystyle 2013x^{4}z+7110x^{5}y^{2}+135x^{7}y-810x^{9}}
\\
\end{pmatrix}$$}

\subsection*{The groups $G_{29}$ and $G_{31}$}
For  the   data  relative  to   $G_{29}$  and  $G_{31}$, see \cite{maschke}.
The  group
$G_{31}$  is  generated  by  the  matrices  $T$  and  $Ue^{2i\pi/8}$  in
Maschke's  notations. 

$G_{29}$  is the  subgroup which  leaves invariant
$\Phi_1$  which  we  take  as   the  first  invariant.  Then  we  choose
$(-1/20736)\det(\Hessian(\Phi_1))=(4F_8-\Phi_1^2)/3$.  We do  not choose
$F_{12}$   but  the   simpler  $((\Phi_1^3-3\Phi_1F_8)/2+F_{12})/108$.
We  do   not  choose   $F_{20}$  but   the  simpler
$(F_{20}-F_8F_{12})/1296$.

For $G_{31}$ we choose $F_8$, $F_{12}$, then as for $G_ {29}$ we choose
$(F_{20}-F_8F_{12})/1296$;  the fourth is still (as in \cite{OT})
$\det(\Hessian(F8))/265531392$.

Basic derivations of $G_{29}$:
{\footnotesize
$$\frac{1}{80}\begin{pmatrix}
320x&640y^{2}&960z+2xy&1600t+8yz\\
640y&4096000t+225280yz+1280xy^{2}&64xz+4x^{2}y&-640tx+16xyz+1536z^{2}\\
960z&-12800tx-640xyz&200t-5yz+3x^{2}z&-10ty+10tx^{2}-8xz^{2}\\
1600t&-51200tz-640txy-1280yz^{2}&-10ty+8tx^{2}-4xz^{2}&72txz-96z^{3}\\
\end{pmatrix}$$}

Basic derivations of $G_{31}$:
$$\frac{1}{270}\begin{pmatrix}
2160x&3240ty&5400z+2xy&6480t-2y^{2}\\
3240y&4860tx^{2}-26244000z^{2}&-9720t+3x^{3}&-11340xz-3x^{2}y\\
5400z&16200xz^{2}-9720t^{2}&-tx-5yz&ty+5x^{2}z\\
6480t&-11340txz-16200yz^{2}&-5ty-2x^{2}z&2xyz+5tx^{2}+5400z^{2}\\
\end{pmatrix}
$$

\subsection*{The group $G_{33}$}
We take the matrices and invariants of \cite[pp 208--209]{bu}
with the corrections indicated in \cite{O}.
The third invariant is taken to be
$\det(\Hessian(J_4))/63700992$ where $J_4$ is the first invariant.
Basic derivations of $G_{33}$:

{\tiny
$$\frac1{128}\begin{pmatrix}
512x&768yz&1280z+(\frac{-4}{3})xy&1536t-4y^{2}&2304u-4ty\\
768y&-663552u+1152x^{2}z&768t-2x^{3}&8064xz-6x^{2}y&-6tx^{2}+23040z^{2}\\
1280z&768tz-1152ux&(\frac{-1}{3})tx+3yz&576u-ty+9x^{2}z&6uy+36xz^{2}-4t^{2}\\
1536t&-3456uy+8064xz^{2}&576u+3ty-5x^{2}z&-15xyz+9tx^{2}+11520z^{2}&-42txz+18u\
x^{2}+108yz^{2}\\
2304u&-3456tu+23040z^{3}&6uy-4xz^{2}&18ux^{2}-12yz^{2}&90uxz-48tz^{2}\\
\end{pmatrix}
$$}

\end{document}